\documentclass{article}
\usepackage[english]{babel}

\usepackage{amsmath,amssymb}
\usepackage[mathscr]{eucal}

\newcommand{\pl}{\partial}
\newcommand{\ol}{\overline}

\newcommand{\HH}{\mathscr{H}}
\newcommand{\WW}{\mathscr{W}}
\newcommand{\Ww}{\mathcal{W}}
\newcommand{\RE}{\mathrm{Re}}
\newcommand{\RR}{\mathbb{R}}
\newcommand{\CC}{\mathbb{C}}
\newcommand{\Sph}{\mathbb{S}}

\begin{document}
\begin{center}
{\Large
Properties of the $l=1$ radial part of the Laplace operator 
in a special scalar product}\\
\vspace{0.3cm}
T.~A.~Bolokhov\\
\vspace{0.3cm}
{\it St.\,Petersburg Department of V.\,A.\,Steklov Mathematical Institute\\
         Russian Academy of Sciences\\
         27 Fontanka, St.\,Petersburg, Russia 191023}
\end{center}

\begin{abstract}
We develop self-adjoint extensions of the $l=1$ radial part of the
Laplace operator in a special scalar product. The product arises as
the transfer of the plain product from 
$ \mathbb{R}^{3} $
into the set of functions parametrizing one of the two components of
the transverse vector field. The similar extensions are treated for
the square of inverse operator of the 
radial part in question.
\end{abstract}


\section*{Introduction}
	The radial part of the Laplace operator
\begin{equation*}
    T_{l} = - \frac{d^{2}}{dr^{2}} + \frac{l(l+1)}{r^{2}} , \quad r\geq 0,
	\quad l = 0,1,\ldots
\end{equation*}
	appears as the
	transfer of the action of the three-dimensional Laplace operator
$ \Delta $
        to the radial components of scalar functions.
        If a function
$ f(\vec{x}) $
        can be represented as a sum of spherical harmonics
$ Y_{lm} $
\begin{equation*}
    f(\vec{x}) = f(\vec{x}(r,\Omega)) = \sum_{0\leq m \leq l}
	\frac{1}{r}f_{lm}(r)
	Y_{lm}(\Omega) , \quad r \geq 0 , \quad \Omega \in \Sph^{2} ,
\end{equation*}
	then the Laplace operator acts on radial components
$ f_{lm}(r) $
    	as the operators
$ T_{l} $ 
\begin{equation*}
    \Delta f(\vec{x})
    = - \sum_{k} \frac{d^{2}}{dx_{k}^{2}} f(\vec{x})
	= \sum_{0\leq m \leq l} \frac{1}{r} T_{l} f_{lm}(r) Y_{lm}(\Omega) .
\end{equation*}
	Similar formulas apply as well for the action of the Laplace operator
	on transverse three-component vector functions, that is, on those
        vector functions that obey the condition
\begin{equation*}
    \vec{f}(\vec{x}) : \quad
	\sum_{k} \frac{\partial f^{k}}{\partial x_{k}} = 0 .
\end{equation*}
	In spherical coordinates these functions can be parameterized
        by two sets of radial components
$ u_{lm}(r) $,
$ \phi_{lm}(r) $
\begin{equation}
\label{Atrexp}
    \vec{f}(\vec{x}) =
        \sum_{1\leq l, |m|\leq l} \bigl( \tilde{l}
	    \frac{u_{lm}}{r^{2}} \vec{\Upsilon}_{lm} +
        \frac{u_{lm}'}{r} \vec{\Psi}_{lm} \bigr)
    +   \sum_{1\leq l, |m|\leq l} \frac{\phi_{lm}}{r} \vec{\Phi}_{lm}  ,
\end{equation}
	where
$ \tilde{l} = \sqrt{l(l+1)} $, and
$ \vec{\Upsilon}(\Omega) $,
$ \vec{\Psi}(\Omega) $,
$ \vec{\Phi}(\Omega) $ --- are (normalized) vector spherical harmonics
\cite{VSH}
\begin{align}
\label{VSH1}
    \vec{\Upsilon}_{lm}(\Omega) = & \frac{\vec{x}}{r} Y_{lm} , \quad
        0 \leq l, \quad |m| \leq l, \\
\label{VSH2}
    \vec{\Psi}_{lm}(\Omega) = & \tilde{l}^{-1} r \vec{\pl} Y_{lm} , \quad
        1 \leq l , \quad |m| \leq l, \\
\label{VSH3}
    \vec{\Phi}_{lm}(\Omega)
	= & \tilde{l}^{-1} (\vec{x} \times \vec{\pl}) Y_{lm},
        \quad 1 \leq l , \quad |m| \leq l .
\end{align}
        In representation 
(\ref{Atrexp})
	the action of the Laplace operator on
$ \vec{f}(\vec{x}) $
        reduces to the action of operators
$ T_{l} $
        on
$ u_{lm}(r) $ and
$ \phi_{lm}(r) $
\begin{equation*}
    \Delta \vec{f}(\vec{x}) =
        \sum_{1\leq l, |m|\leq l} \bigl( \tilde{l}
	    \frac{T_{l}u_{lm}}{r^{2}} \vec{\Upsilon}_{lm} +
        \frac{(T_{l}u_{lm})'}{r} \vec{\Psi}_{lm} \bigr)
    +   \sum_{1\leq l, |m|\leq l} \frac{T_{l}\phi_{lm}}{r} \vec{\Phi}_{lm}  ,
\end{equation*}
	for details see \cite{Lapl}.

        The plane scalar product in the space of functions on
$ \RR^{3} $
\begin{equation}
\label{fgprod}
    (\vec{f},\vec{g}) = \int_{\RR^{3}} \ol{\vec{f}(\vec{x})} \cdot
	\vec{g}(\vec{x})\, d^{3}x
    = \int_{\RR^{3}} \sum_{j} \ol{f^{j}(\vec{x})} g^{j}(\vec{x}) \,d^{3}x
\end{equation}
	naturally transfers onto the sets of functions
$ u_{lm}(r) $, 
$ \phi_{lm}(r) $
\begin{align}
\label{prod1}
    \bigl(\frac{\phi_{lm}}{r}\vec{\Phi}_{lm}, &
	\frac{\tilde{\phi}_{l'm'}}{r} \vec{\Phi}_{l'm'}\bigr) =
	\delta_{ll'} \delta_{mm'} \,
    \int_{0}^{\infty} \ol{\phi_{lm}} \tilde{\phi}_{l'm'}\,dr 
    \delta_{ll'}\delta_{mm'} \,
    \equiv (\phi_{lm}, \tilde{\phi}_{l'm'}) ,\\
\nonumber
    \bigl(\tilde{l}\frac{u_{lm}}{r^{2}}\vec{\Upsilon}_{lm} &
	    + \frac{u'_{lm}}{r} \vec{\Psi}_{lm} ,
	\tilde{l'}\frac{v_{l'm'}}{r^{2}}\vec{\Upsilon}_{l'm'}
	    + \frac{v'_{l'm'}}{r} \vec{\Psi}_{l'm'} \bigr) = \\
\label{prod2}
    & = \delta_{ll'} \delta_{mm'} \,
	\int_{0}^{\infty} \bigl(\ol{u'}v' + \frac{l(l+1)}{r^{2}} \ol{u}v
	\bigr) dr
    \equiv \delta_{ll'} \delta_{mm'} \, \langle u, v \rangle_{l} .
\end{align}
	It follows immediately that the properties of the transverse Laplace operator
        in spherical coordinates are determined by the properties of radial operators
$ T_{l} $
        in scalar products
(\ref{prod1}) and
(\ref{prod2}).
        The first case matches the known example of the action of the
	scalar Laplacian in the subspaces corresponding to 
$ l\geq1 $,
        while the second case turns out to be not as trivial.
        In reference
\cite{Lapl}
	the quadratic form of the transverse Laplace operator, 
	transferred to functions
$ u_{lm} $
\begin{equation}
\label{lmform}
    Q_{\Delta}(\vec{f}) =
    \bigl(\tilde{l}\frac{u_{lm}}{r^{2}}\vec{\Upsilon}_{lm} 
	    + \frac{u'_{lm}}{r} \vec{\Psi}_{lm} , \,
	\tilde{l'}\frac{T_{l} u_{l'm'}}{r^{2}}\vec{\Upsilon}_{l'm'}
	    + \frac{(T_{l} u_{l'm'})'}{r} \vec{\Psi}_{l'm'} \bigr) ,
\end{equation}
	is studied as a form in the scalar product
\begin{equation}
\label{plainprod}
    (u,v) = \int_{0}^{\infty} \ol{u(r)} v(r) dr .
\end{equation}
	The corresponding generating operator
	appears to be a fourth order differential operator.
	For
$ l=1,2 $
        such operators are symmetric with deficiency indices 
$ (1,1) $
	and have nontrivial self-adjoint extensions.
        It turns out, however, that for
$ l=1 $
        the quadratic form
(\ref{lmform})
        and operator
$ T_{1} $
        have nontrivial extensions also in the scalar product
        transferred from
$ \RR^{3} $
\begin{equation}
\label{angleprod}
    \langle u, v \rangle_{l} =
    \int_{0}^{\infty} \bigl(\ol{u'}v' + \frac{l(l+1)}{r^{2}} \ol{u}v
	\bigr) dr .
\end{equation}
	(This idea was proposed in \cite{Lapl}, but was not developed therein.)
	In particular, such an assertion allows one to study a physical model
        analogous to 
\cite{BF}
	in the case of a transverse field.
        The description of the resolvent and spectral properties of
    extended operator
$ T_{1} $
        is given in the second section of the present work.

        In the third section we study extensions of the operator corresponding
        to the quadratic form
\begin{equation}
\label{t1sform}
    Q^{-1}(u) =
    \iint_{0}^{\infty} \ol{u(r)} T_{l}^{-1}(r,s) u(s) \,dr\,ds , \quad l=1 ,
\end{equation}
	considered in scalar product
(\ref{angleprod}).
	In scalar product
(\ref{plainprod})
	this form is a quadradic form of the operator which is
	inverse to the essentially self-adjoint operator
$ T_{l}|_{l=1} $,
	and so it does not have extensions.
	This form does not possess extensions in scalar product
(\ref{angleprod}) also.
	However, in the latter case the form is an extension of another
	form which is defined on some set, dense in the norm
(\ref{angleprod}).

        Quadratic forms of the kind 
(\ref{lmform}) and
(\ref{t1sform})
	appear, correspondingly, in the potential and kinetic parts of the
        Hamiltonian for the electromagnetic field in the Coloumb gauge, when written
	in spherical coordinates.
	Therefore explicit expressions for the corresponding self-adjoint operators
        are beneficial for a more profound understanding of the quantization
        of the latter model.
    Although we should note, that the scalar product
(\ref{angleprod})
    is unphysical for the kinetic part of the Hamiltonian. Thus
    the third section
    is provided as an example of resolvent technique for the
    inverse operators.

\section{Notations and preliminaries}
\label{notations}
	For simplicity we will omit the orbital momentum index
	in the notations for operators and the scalar product henceforth.
	Furthermore, under the notation of an ``operator
$ T $''
	we will assume the differential operation
\begin{equation*}
    T = -\frac{d^{2}}{dr^{2}} + \frac{l(l+1)}{r^{2}} \Bigr|_{l=1}
	= -\frac{d^{2}}{dr^{2}} + \frac{2}{r^{2}} 
\end{equation*}
	without specification of the domain.
	The parentheses 
(\ref{plainprod}) will designate the corresponding integral --- 
	just on the set	where it converges.

        In order to facilitate construction of the 
$ T $-invariant functions (``eigenfunctions'' in a broad sence)
	let us introduce covariant derivatives
\begin{equation*}
    D \equiv r \frac{d}{dr} r^{-1} = \frac{d}{dr} -\frac{1}{r} , \quad
    D^{*} \equiv - r^{-1}\frac{d}{dr} r = - \frac{d}{dr} - \frac{1}{r} .
\end{equation*}
	It is not difficult to see that
$ T $, $ D $ and $ D^{*} $
	obey the following algebraic identities
\begin{equation}
\label{TDrel}
    T = DD^{*} , \quad TD = D D^{*} D , \quad D^{*}D = - \frac{d^{2}}{dr^{2}}.
\end{equation}
	Below we will be using the following expansion and integrals which
	involve derivative
$ D $
	acting on an exponent:
\begin{align}
\label{Pf1}
    \sum_{n} \alpha_{n} D e^{\sigma_{n} r} \sim
	\sum_{n} \alpha_{n} \bigl( - r^{-1}
	  &+ \sigma_{n}^{2} \frac{r}{2}
	+ \sigma_{n}^{3} \frac{r^{2}}{3} + \ldots \bigr),
	\quad r\to 0\\
\label{Pf2}
    \int_{0}^{\infty} \sum_{n} \frac{\alpha_{n}}{r} De^{\sigma_{n}r} \, dr & =
	- \sum_{n} \alpha_{n} \sigma_{n} ,\quad \sum_{n} \alpha_{n} = 0 , \\
\label{Pf3}
    \int_{0}^{\infty} \sum_{n} r^{2} \alpha_{n} De^{\sigma_{n}r} \, dr & =
- 3 \sum_{n} \frac{\alpha_{n}}{\sigma_{n}^{2}} .
\end{align}
	We also introduce the notation
\begin{equation*}
    T^{-1} = T^{-1}(r,s) =
	\frac{r^{2}}{3s} \theta(s-r) + \frac{s^{2}}{3r} \theta(r-s) 
\end{equation*}
	for the kernel of the operator that is inverse to the
	essentially self-adjoint operator 
$ T $
	in scalar product
(\ref{plainprod}).
	This property of the kernel of
$ T^{-1} $
	is confirmed by a direct calculation
\begin{equation*}
    \bigl(-\frac{d^{2}}{dr^{2}} + \frac{2}{r^{2}} \bigr) \bigl(
	\frac{r^{2}}{3s} \theta(s-r) + \frac{s^{2}}{3r} \theta(r-s)  
    \bigr) = \delta(r-s) .
\end{equation*}
	We also denote as
$ \HH $
	the Hilbert space of functions on the half-axis which is
        generated by scalar product
$ \langle \cdot , \cdot \rangle $
\begin{equation*}
    \HH = \{u(r): \int_{0}^{\infty} \bigr(|u'|^{2}
	+\frac{2}{r^{2}}|u|^{2}\bigr) dr < \infty \} .
\end{equation*}
	This space is comprised of absolutely continuous functions
	vanishing at zero, for which the above integral is finite.
	If a function
$ v $
	from space
$ \HH $
	is twice differentiable then one can write a formal equality
        for product
$ \langle u, v \rangle $
\begin{equation}
\label{apc}
    \langle u,v\rangle = \int_{0}^{\infty}
	\bigl(\ol{\frac{du}{dr}}\frac{dv}{dr}+\frac{2}{r^{2}}\ol{u}v\bigr)dr
    = \int_{0}^{\infty} \ol{u}\bigl(-\frac{d^{2}v}{dr^{2}}+\frac{2v}{r^{2}}
	\bigr) dr = (u,Tv) .
\end{equation}
	As we have pointed out above, the expression on the right hand side
	of this equality in general is not a scalar product, as either of
        functions
$ u $	and
$ Tv $
	may not be square integrable.

\subsection{Expressions for the resolvent}
	We shall use some consequences of the Stone's formula
	for the resolvent of a self-adjoint operator (see, {\it e.g.}, 
\cite{Resolvent}).
	Let
$ R_{A}(r,s;w) $ be the kernel of the resolvent of a self-adjoint operator
$ A $
	with a positive continuous spectrum, and acting on functions
        of a real variable:
\begin{equation*}
    (A - w) R_{A}(r,s;w) = \delta(r-s) , \quad w \in \CC \setminus \sigma(A) .
\end{equation*}
	Using this kernel one can write an identity decomposition (the completeness relation)
\begin{equation}
\label{Acompl}
    \delta(r-s) = \frac{1}{2\pi i} \int_{0}^{\infty} \bigl(
	R_{A}(\eta)-R_{A}(e^{2\pi i}\eta) \bigr) \, d\eta
    +\sum_{n} \lim_{w\to w_{n}} R_{A}(w) (w_{n}-w) ,
\end{equation}
	where
$ w_{n} $ are the poles of the resolvent
$ R(w) $.
	If the spectrum of
$ A $
	is of multiplicity one
	then the integrand in this formula can be written as
        a square of real-valued ``eigenfunctions of the continuous spectrum''
\begin{equation}
\label{rescs}
    R_{A}(r,s;\eta) - R_{A}(r,s;e^{2\pi i}\eta)
	= 2\pi i p_{\eta}(r) p_{\eta}(s) ,
\end{equation}
	and the sum of residues --- as a sum of squares of discrete spectrum
	eigenfunctions
\begin{equation}
\label{resps}
    \sum_{n} \lim_{w\to w_{n}} R(r,s;w) (w_{n}-w) = \sum_{n} q_{n}(r) q_{n}(s) .
\end{equation}
	Besides the completeness relations
(\ref{Acompl}),
	the following orthogonality relations hold for
$ q_{n}(r) $,
$ p_{\eta}(r) $,
\begin{gather*}
    \int p_{\eta}(r) p_{\xi}(r) dr = \delta(\eta -\xi) ,\\
    \int q_{n}(r) p_{\eta}(r) dr = 0 , \quad
    \int q_{n}(r) q_{m}(r) dr = \delta_{nm} .
\end{gather*}
	In the case when operator
$ A $
	and resolvent
$ R(w) $
	are self-adjoint in a space with a scalar product
$ \langle \cdot , \cdot \rangle $
	obeying
(\ref{apc}),
	expansions
(\ref{rescs}), (\ref{resps})
	get an additional factor of operator
$ T $
	acting on the ``right-hand side'' variable
$ s $:
\begin{align*}
    R(r,s;\eta) - R(r,s;e^{2\pi i}\eta) =& 2\pi i p_{\eta}(r) Tp_{\eta}(s) , \\
    \lim_{w\to w_{n}} R(r,s;w) (w_{n}-w) =& q_{n}(r) Tq_{n}(s) ,
\end{align*}
	while in the orthogonality relations the integral changes into the bracket
$ \langle \cdot , \cdot \rangle $.

\section{Operator $ T $ in scalar product
$ \langle \cdot , \cdot \rangle $}
	Let us take the space of smooth functions
$ \Ww_{0} $,
\begin{equation*}
    \Ww_{0} = \{u(r): \; u \in \HH, \: Tu \in \HH,\: u''(0)=0 \}
\end{equation*}
	as the initial domain of the symmetric operator.
	We remark that for
$ Tu $  to belong to
$ \HH $,
$ u $ being a smooth function,
	only the first derivative
$ u'(0) = 0 $
        needs to vanish at zero,
	while, because of the relation
 $ T r^{2} = 0 $,
	the second one needs not.
	For that reason the condition
$ u''(0) = 0 $
	in
$ \Ww_{0} $
	is demanded separately.

	Let us define a symmetric operator
$ \tilde{T} $
	as a restriction of operator
$ T $
	onto the set
$ \Ww_{0} $:
$ \tilde{T} = T|_{\Ww_{0}} $.
	Its symmetricity is established by integrating by parts
        and expanding the smooth functions in a series in
        a vicinity of zero:
\begin{equation*}
    \langle u, Tv \rangle = \langle Tu, v \rangle + 4(u'v''-u''v')|_{r=0}
	= \langle Tu, v \rangle .
\end{equation*}

\subsection{Deficiency indices}
	In order to calculate the deficiency indices of operator
$ \tilde{T} $
	it is necessary to study the kernels of the adjoint operators
$ \tilde{T}^{*} \mp i\rho^{2} $,
$ \rho > 0 $.
	Let
$ c_{\pm} $
	be any two vectors from the corresponding kernels;
	then for any
$ v \in \Ww_{0} $
	the following equality holds
\begin{equation}
\label{adjc}
    \langle c_{\pm} , (\tilde{T}\pm i\rho^{2}) v \rangle = 0, \quad
	v \in \Ww_{0} .
\end{equation}
	It is not hard to see that as
$ c_{\pm} $
    one can take the vectors
\begin{equation*}
    c_{\pm} = D \exp \{e^{\mp \frac{3\pi i}{4}}\rho r\}+ r^{-1} 
\end{equation*}
	--- these functions belong to 
$ \Ww_{0} $
	and fulfill relation
(\ref{adjc}):
\begin{align*}
    \langle D& \exp \{e^{\mp \frac{3\pi i}{4}}\rho r\}+ r^{-1} ,
	(\tilde{T}\pm i\rho^{2})g \rangle = 
    \bigl(T (D\exp \{e^{\frac{3\pi i}{4}}\rho r\}+r^{-1}),
	(T\pm i\rho^{2})g \bigr) \\
    &= \bigl( (T\pm i\rho^{2}) T
	(D\exp \{e^{\frac{3\pi i}{4}}\rho r\}+r^{-1}), g \bigr) 
    = \bigl( (T\pm i\rho^{2}) T
	D\exp \{e^{\frac{3\pi i}{4}}\rho r\}, g \bigr) \\
    &= \bigl( T (T\pm i\rho^{2})
	D\exp \{e^{\frac{3\pi i}{4}}\rho r\}, g \bigr) 
    = \bigl( T D(-\frac{d^{2}}{dr^{2}} \pm i\rho^{2})
	\exp \{e^{\frac{3\pi i}{4}}\rho r\}, g \bigr) = 0.
\end{align*}
	Here the first equality is a re-writing of relation
(\ref{apc}),
	while the second one is found by integrating by parts
        and demands that
$ g''(0) = 0 $.
	The third equality follows from the fact that
$ T r^{-1} = 0 $,
	the fourth one is merely a commutation of the operations
$ T $ and
$ T\pm i\rho^{2} $,
	while the fifth one is a consequence of relation
(\ref{TDrel}).

	The absence of other elements in the kernels of the operators
$ \tilde{T}^{*} \mp i\rho^{2} $
	follows from the theory of differential equations for distributions
(see, {\it e.g.}
\cite{RS1}).
	The corresponding fourth order equation has four
	linearly independent solutions.
	One of them is exhibited above, while the other three,
$ r^{-1} $,
$ r^{2} $,
$ D\exp\{e^{\pm\frac{\pi i}{4}}\rho r\} $,
	are either divergent at zero, or not integrable at infinity,
	even in the norm
$ \langle \cdot , \cdot \rangle $.

	We therefore conclude that the symmetric operator
$ \tilde{T} $
	has deficiency indices
$ (1,1) $
	and, consequently, possesses non-trivial self-adjoint extensions.

\subsection{Self-adjoint extensions of operator
$ \tilde{T} $}
	The standard way to construct self-adjoint extensions of a symmetric operator
	is the Cayley transform
(see, {\it e.g.} \cite{RS}).
	We shall not dwell upon this technique here but, instead, provide
	the domains and the extended operators.
	Explicitly, let
\begin{equation}
\label{bcond}
    \Ww_{\kappa} = \{u(r): \; u \in \HH, \: T_{\kappa}u \in \HH,\:
	3u''(0) = 4\kappa u'(0) \}, \quad \kappa \in \RR
\end{equation}
	be the domain where the extended operator
$ T_{\kappa} $
	acts as a differential operation
$ T $
	with a non-local addition
\begin{equation}
\label{Tk}
    T_{\kappa} u = Tu - \frac{2}{r} u'(0)
    = -\frac{d^{2}u}{dr^{2}} + \frac{2}{r^{2}} u -\frac{2}{r}u'(0) .
\end{equation}
	It is obvious that
$ \tilde{T} \subset T_{\kappa} $,
	since
$ \Ww_{0} \subset \Ww_{\kappa} $ and
$ T_{\kappa} u = \tilde{T} u $, whenever
$ u \in \Ww_{0} $.
	The symmetricity of operator
$ T_{\kappa} $ acting on
$ \Ww_{\kappa} $
    can be ensured via an explicit calculation, as a result of which we find,
\begin{equation*}
    \langle u, T_{\kappa}v \rangle = \langle T_{\kappa} u, v \rangle
	+ 5(u'v''-u''v')|_{r=0}	= \langle T_{\kappa}u, v \rangle .
\end{equation*}
	Here
$ u'v''-u''v'|_{r=0}=0 $
	if
$ u $ and $ v $
	both obey the boundary conditions from
$ \Ww_{\kappa} $.
	Self-adjointness of 
$ T_{\kappa} $ on
$ \Ww_{\kappa} $
	follows from the expansion of 
$ c_{\pm} $
	near zero which can be obtained using Eq.~\eqref{Pf1}:
\begin{equation*}
    c_{\pm}(r) \sim \mp i\rho^{2} \bigl(\frac{r}{2}
	+ e^{\mp\frac{3\pi i}{4}}\rho \frac{r^{2}}{3} +\ldots \bigr)
\end{equation*}
	--- such functions do not satisfy the boundary conditions from
$ \Ww_{\kappa} $
	for any
$ \rho $ and
$ \kappa $,
	and, therefore, the kernels of the adjoint operators
$ T_{\kappa}^{*} \mp i\rho^{2} $,
	which are embedded in the kernels of
$ \tilde{T}^{*} \mp i\rho^{2} $,
	are empty.

\subsection{Resolvent of operator $ T_{\kappa} $}
	One way to study the spectral properties of 
$ T_{\kappa} $
	is to exploit the properties of its resolvent kernel.
        Let us seek this kernel as a function
$ R(r,s;z) $
        of two positive and one complex variables, satisfying the equation
\begin{equation}
\label{Treseq}
    (T_{\kappa}-z^{2}) R(r,s;z) = \delta(r-s) ,\quad 0 <\arg z < \pi ,
\end{equation}
	and the boundary conditions
$ \Ww_{\kappa} $
	in
$ r $.
        We also require that this function should obey a complex symmetricity
        condition with respect to scalar product
(\ref{angleprod}),
        and also that it should fall off at infinity in variables
$ r $ and
$ s $.
        Such a function can be constructed by means of the two solutions
$ h(r) $,
$ g(r) $
        of a homogeneous equation
\begin{gather*}
    (T-z^{2}) h(r) = 0, \quad (T-z^{2}) g(r) = 0\\
    h(r) = D e^{-izr} + \beta(z) D e^{izr} , \quad
    g(r) = D e^{izr} .
\end{gather*}
	Explicitly,
\begin{equation}
\label{resT}
    R(r,s;z) = \frac{1}{W}\bigl(h(r)g(s)\theta(s-r) + h(s)g(r) \theta(r-s)
	+\frac{1+\beta}{r} g(s) \bigr) ,
\end{equation}
	where
\begin{equation*}
    W(z) = h'(r) g(r) - h(r)g'(r) = -2iz^{3} .
\end{equation*}
	Let us begin by checking that expression
(\ref{resT})
        satisfies the resolvent boundary conditions.
        The exponents
$ g(r) $, $ g(s) $
        which are selected by the 
$ \theta $-function
	at
$ r\to\infty $ or at
$ s\to\infty $
        and function
$ r^{-1} $
        fall off at infinity if
$ 0 < \arg z < \pi $.
        At a given
$ s $,
        the second term in
(\ref{resT})
        vanishes in a vicinity of zero in
$ r $,
        while the rest of the terms have the following expansion
\begin{equation}
\label{boundr}
    \frac{g(s)}{W} \bigl(h(r)+\frac{1+\beta}{r}\bigr)
	\sim \frac{g(s)}{W} \bigl(
	-\frac{z^{2}}{2}(1+\beta)r + \frac{iz^{3}}{3}(1-\beta)r^{2} + \ldots
    \bigr),
\end{equation}
	which can be seen to satisfy the boundary conditions from
$ \Ww_{\kappa} $,
        provided that
\begin{equation*}
    \beta(z) = \frac{z-i\kappa}{z+i\kappa} .
\end{equation*}
	Expansion
(\ref{boundr})
	also shows that the last (non-local) term in the operator
$ T_{\kappa} $
        when applied to
$ R(r,s;z) $
        gives
\begin{equation}
\label{addt}
    - \frac{2}{r} \frac{\partial}{\partial \tilde{r}}
      R(\tilde{r},s;z)
	\bigr|_{\tilde{r}=0} = z^{2} \frac{(1+\beta)}{Wr} g(s).
\end{equation}
	Next we note that the first two terms in
(\ref{resT})
        form a standard combination of two solutions of a homogeneous equation
        for a resolvent of a second order differential operator,
        and therefore they satisfy the functional equation
\begin{equation}
\label{resc}
    (T-z^{2}) \frac{1}{W}\bigl(h(r)g(s)\theta(s-r) + h(s)g(r) \theta(r-s)
	\bigr) = \delta(r-s).
\end{equation}
	The last term in
(\ref{resT}) is annihilated by
$ T $,
        and we get
\begin{equation}
\label{addr}
    (T-z^{2})\frac{1+\beta}{Wr} g(s) = -z^{2} \frac{1+\beta}{Wr} g(s) ,
\end{equation}
	which cancels the contribution
(\ref{addt})
        of the non-local term in
$ T_{\kappa} $.
	Altogether, summing up terms
(\ref{addt}),
(\ref{resc}) and
(\ref{addr})
	we find that expression
(\ref{resT})
	does satisfy equation
(\ref{Treseq}).

        To check the complex symmetricity condition of the kernel
(\ref{resT})
        with respect to scalar product
(\ref{angleprod}),
	which in terms of variable
$ z $ 
	looks as
\begin{equation*}
    R^{*}(r,s;z) = R(s,r;-\ol{z}),
\end{equation*}
        let us represent
$ R(r,s;z) $
        as a result of the action of
$ T^{-1} $
        on a certain kernel
$ S(r,s;z) $
\begin{align*}
    R(r,s;z) &= \int T^{-1}(r,t)\Bigl( \frac{i}{2z}\bigl(
	h(t)g(s)\theta(s-t) + h(s)g(t)\theta(t-s)\bigr) \\
    & +\delta(t-s)\Bigr) \,dt = T^{-1} S ,
	\quad \ol{S(s,r;z)} = S(r,s;-\ol{z})
\end{align*}
	which is complex symmetric by construction.
	Then, given arbitrary
$ u,v \in \Ww_{\kappa} $
        we can use relation
(\ref{apc})
        and write
\begin{align*}
    \langle u,R(z)v \rangle &= (u,TR(z)v) = (u,S(z)v) \\
	&= (S(-\ol{z})u,v) 
	= (TR(-\ol{z})u,v) = \langle R(-\ol{z})u, v \rangle ,
\end{align*}
	where the parentheses denote the integral (\ref{plainprod}),
	as we specified in Section~\ref{notations}.

\subsection{Spectral decomposition of
$ T_{\kappa} $}
	In the case when the resolvent
$ R $
	is given in terms of variable
$ z = \sqrt{w} $
	and satisfies equation
(\ref{Treseq}),
	identity decomposition
(\ref{Acompl})
	can be cast into the following form,
\begin{align}
\label{Tcompl}
    &\delta(r-s) = \frac{1}{2\pi i} \int_{0}^{\infty} \bigl(
	R(r,s;\lambda)-R(r,s;-\lambda) \bigr) 2\lambda \, d\lambda \\
\nonumber
    &+ \sum_{n} 2z_{n} \lim_{z\to z_{n}} R(r,s;z) (z_{n}-z)
= \int_{0}^{\infty} \tilde{p}_{\lambda}(r) T\tilde{p}_{\lambda}(s) \,d\lambda
	+ \sum_{n} \tilde{q}_{n}(r) T\tilde{q}_{n}(s) .
\end{align}
	Notably, the orthogonality conditions for 
$ \tilde{p}_{\lambda} $
	are valid in terms of parameter
$ \lambda $ as well,
\begin{gather}
\label{ortrel1}
    \int \tilde{p}_{\lambda}(r) T \tilde{p}_{\mu}(r) dr
	= \delta(\lambda -\mu) ,\\
\label{ortrel2}
    \int \tilde{q}_{n}(r) T \tilde{p}_{\lambda}(r) dr = 0 , \quad
    \int \tilde{q}_{n}(r) T \tilde{q}_{m}(r) dr = \delta_{nm} .
\end{gather}
	To calculate the eigenfunctions of the discrete spectrum we examine
	the poles of the resolvent
$ R(\lambda) $: its only pole resides in the coefficient
$ \beta(z) $
	and is located at
\begin{equation*}
    z_{0} = -i\kappa.
\end{equation*}
	It falls into the upper half-plane only when
$ \kappa < 0 $.
	In that case,
\begin{align*}
    \lim_{z\to z_{0}} 2z_{0} R(r,s;z) (z_{0}-z)=&\frac{4z_{0}^{2}}{2iz_{0}^{3}}
	\bigl(De^{\kappa r}
	    + \frac{1}{r}\bigr) D e^{\kappa s} \\
	=& -\frac{2}{\kappa^{3}}
    \bigl(De^{\kappa r} + \frac{1}{r}\bigr) T
	\bigl(D e^{\kappa s} + \frac{1}{s}\bigr) ,
\end{align*}
	and thus we obtain a simple expression for the
	normalized discrete spectrum eigenfunction,
\begin{equation}
\label{Tqr}
    \tilde{q}(r) = \tilde{q}_{\kappa}(r)
    = \sqrt{-\frac{2}{\kappa^{3}}}
	\bigl(D e^{\kappa r} + \frac{1}{r}\bigr) .
\end{equation}
	To calculate the spectral density of the continuous component
    of the spectrum
	we perform the following derivation,
\begin{align*}
    \frac{2\lambda}{2\pi i} &\bigl(
	R(r,s;\lambda)-R(r,s;-\lambda) \bigr)
    = \frac{-2\lambda}{2\pi i 2i\lambda^{3}} \bigl(
	\frac{1+\beta_{\lambda}}{r} De^{i\lambda s} 
	+\frac{1+\beta_{\lambda}^{-1}}{r} De^{-i\lambda s} \\
&+(De^{-i\lambda r}+\beta_{\lambda}De^{i\lambda r}) De^{i\lambda s} \theta(s-r)
+(De^{-i\lambda s}+\beta_{\lambda}De^{i\lambda s}) De^{i\lambda r}
    \theta(r-s)\\
&+(De^{i\lambda r}+\beta_{\lambda}^{-1} De^{-i\lambda r}) De^{-i\lambda s}
    \theta(s-r)
+(De^{i\lambda s}+\beta_{\lambda}^{-1} De^{-i\lambda s}) De^{-i\lambda r}
    \theta(r-s) \bigr)\\
=& \frac{1}{2\pi \lambda^{2}} \bigl(
	\frac{1+\beta_{\lambda}}{r} De^{i\lambda s} 
	+\frac{1+\beta_{\lambda}^{-1}}{r} De^{-i\lambda s} \\
&+ De^{-i\lambda r} De^{i\lambda s} + De^{i\lambda r} De^{-i\lambda s}
    + \beta_{\lambda} De^{i\lambda r} De^{i\lambda s}
    + \beta_{\lambda}^{-1} De^{-i\lambda r} De^{-i\lambda s} \bigr) \\
=& \frac{1}{2\pi \lambda^{2}} \bigl(
    \beta_{\lambda}^{\frac{1}{2}} De^{i\lambda r}
    + \beta_{\lambda}^{-\frac{1}{2}} De^{-i\lambda r}
+ \frac{\beta_{\lambda}^{\frac{1}{2}}+\beta_{\lambda}^{-\frac{1}{2}}}{r}\bigr)
    \bigl( \beta_{\lambda}^{\frac{1}{2}} De^{i\lambda s}
    + \beta_{\lambda}^{-\frac{1}{2}} De^{-i\lambda s} \bigr) ,
\end{align*}
	where
\begin{equation*}
    \beta_{\lambda} =\beta(\lambda)
	= \frac{\lambda-i\kappa}{\lambda +i\kappa} \equiv e^{2i\zeta} .
\end{equation*}
	This spectral density produces the following real-valued
	``continuous spectrum eigenfunction''
$ \tilde{p}_{\lambda}(r) $
	from decomposition
(\ref{Tcompl}),
\begin{equation}
\label{Tpl}
    \tilde{p}_{\lambda}(r) = \tilde{p}_{\lambda,\kappa}(r)
	= \frac{1}{\sqrt{2\pi}\lambda^{2}} \bigl(
	De^{i\zeta +i\lambda r} + De^{-i\zeta -i\lambda r}
	    +\frac{e^{i\zeta}+e^{-i\zeta}}{r}\bigr) .
\end{equation}

\subsection{The quadratic form}
    Concluding this section we would like to derive the expression
    for the pullback of the forms
$ \langle u, T_{\kappa} u \rangle $ 
    into the three-dimensional space. We shall assume here that repeated
    indices
$ k $, 
$ j $, $ l $ and $ m $
    are summed.
    Let us start with the quadratic form
(\ref{lmform})
    for the transverse field
$ \vec{f} $
    parametrized as in
(\ref{Atrexp})
\begin{equation}
\label{Qform}
    Q_{\Delta}(\vec{f}) = -\int_{\RR^{3}} f^{j}(\vec{x})
\frac{\partial^{2} f^{j}(\vec{x})}{\partial x_{k}^{2}}\, d^{3} x =
(\phi_{lm},T_{l}\phi_{lm})
    +\langle u_{lm}, T_{l}u_{lm}\rangle_{l} 
\end{equation}
    and change the forms
$ \langle u_{1m}, T_{1} u_{1m} \rangle_{1} $,
    corresponding to
$ l=1 $,
    to their extensions
$    \langle u_{1m}, T_{1\kappa}u_{1m} \rangle_{1} $,
    where
\begin{equation*}
    T_{1\kappa}u_{1m}(r) = T_{1}u_{1m}(r) -2r^{-1}u'_{1m}(0).
\end{equation*}
    The scalar products in the extended form can be written as the common limit
\begin{align*}
Q_{\Delta}^{\kappa}&(\vec{f}) 
= (\phi_{lm},T_{l}\phi_{lm})
+\langle u_{lm}, T_{l}u_{lm}\rangle_{l}\bigr|_{2\leq l}
+\langle u_{1m}, T_{1\kappa}u_{1m}\rangle_{1} = \\
&= \lim_{\rho\to 0} \int_{\rho}^{\infty} \Bigl(\phi_{lm} T_{l}\phi_{lm} 
	+ \bigl(
    \frac{du_{lm}}{dr} \frac{d}{dr} T_{l} u_{lm}
	+ \frac{l(l+1)}{r^{2}} u_{lm}T_{l}u_{lm} \bigr)_{2\leq l} +\\
	&+ \bigl(
    \frac{du_{1m}}{dr} \frac{d}{dr} T_{1} u_{1m}
	+ \frac{2}{r^{2}} u_{1m} T_{1} u_{1m} \bigr) 
    + 2\bigl(\frac{u'_{1m}(r)}{r^{2}}-\frac{2u_{1m}(r)}{r^{3}}\bigr)u'_{1m}(0)
    \Bigr)dr.
\end{align*}
    Everything containing
$ T_{l} $
    can be transferred back to the function
$ \vec{f}_{m} $
    by the equation 
(\ref{Qform}),
    restricted to the complement of the ball
$ B_{\rho} $,
    and then integrated by parts
\begin{align}
\label{qftrans}
Q_{\Delta}^{\kappa}(\vec{f}) 
    &= - \int_{\RR^{3}\setminus B_{\rho}} f^{j}(\vec{x})
    \frac{\partial^{2} f^{j}(\vec{x})}{\partial x_{k}^{2}} \,d^{3}x
    - 2u'_{1m}(0) \frac{u_{1m}(r)}{r^{2}}\bigr|_{r=\rho} = \\
\nonumber
    =& \int_{\RR^{3}\setminus B_{\rho}}
	\bigl( \frac{\partial f^{j}(\vec{x})}{\partial x_{k}}
	\bigr)^{2} d^{3}x
    + \int_{\partial B_{\rho}}
      \partial_{k} f^{j}(\vec{x}) f^{j}(\vec{x}) \,d^{2}\sigma_{k}
    - 2u'_{1m}(0) \frac{u_{1m}(\rho)}{\rho^{2}} ,
\end{align}
    here
$ d^{2}\vec{\sigma} $
    is the vector from the divergence theorem, normal to the element 
$ d^{2}\sigma $
    of the boundary
$ \partial B_{\rho} $.
    Substituting
$ \vec{x} = \vec{x}(\rho,\Omega) $
    and
$ d^{2}\vec{\sigma} = \vec{x} \rho d^{2}\Omega $,
    where
$ \Omega $
    is a point on the unit sphere
$ \Sph^{2} $,
    the integral in the second term is written as
\begin{equation}
\label{sterm}
    \int_{\partial B_{\rho}}
	\partial_{k} f^{j}(\vec{x}) f^{j}(\vec{x}) \,d^{2}\sigma^{k}
    = \int_{\Sph^{2}}
    \partial_{k} f^{j}(\vec{x})f^{j}(\vec{x}) x_{k} \rho\,d^{2}\Omega .
\end{equation}
    The partial derivative
$ \partial_{k}f^{j} $,
    using the definitions
(\ref{VSH1}), (\ref{VSH2}),
    can be represented by 5 terms
\begin{multline*}
    \partial_{k}f^{j}(\vec{x}(\rho,\Omega))
	= \frac{\tilde{l}x_{k}}{\rho^{3}}
	\bigl(\rho u_{lm}'(\rho)-3u_{lm}(\rho)\bigr)
    \Upsilon^{j}_{lm}(\Omega)
	+ \frac{\tilde{l}u_{lm}(\rho)}{\rho^{3}} \delta_{kj} Y_{lm}(\Omega) +\\
    + \frac{\tilde{l}^{2} u_{lm}(\rho)}{\rho^{4}} x^{j} \Psi^{k}_{lm}(\Omega)
    + \frac{x_{k}}{\rho^{2}}u_{lm}''(\rho) \Psi^{j}_{lm}(\Omega)
    + \tilde{l}^{-1}u_{lm}'(\rho) \partial_{k}\partial_{j} Y_{lm}(\Omega) .
\end{multline*}
    It is not hard to check, that the following relations hold
\begin{align*}
    x_{k}& \Psi^{k}_{lm}(\Omega) = x_{k} \rho \partial_{k} Y_{lm}(\Omega)
	= 0, \quad \vec{x} = \vec{x}(\rho,\Omega),\\
    \int_{\Sph^{2}} & \Upsilon^{j}_{lm}(\Omega) x_{k}
	\partial_{j}\partial_{k} Y_{l'm'}(\Omega)\, d^{2}\Omega  
    =- \int_{\Sph^{2}} \Upsilon^{j}_{lm}(\Omega) \delta^{kj}
	\partial_{k} Y_{l'm'}(\Omega)\, d^{2}\Omega =\\
    &= - \int_{\Sph^{2}} \rho^{-1} \Upsilon^{j}_{lm}(\Omega) 
	\Psi^{j}_{l'm'}(\Omega)\, d^{2}\Omega = 0 \\
    \int_{\Sph^{2}} & \Psi^{j}_{lm}(\Omega) x_{k}
	\partial_{j}\partial_{k} Y_{l'm'}(\Omega)\, d^{2}\Omega 
    =- \int_{\Sph^{2}} \Psi^{j}_{lm}(\Omega) \delta^{kj}
	\partial_{k} Y_{l'm'}(\Omega)\, d^{2}\Omega =\\
    &= - \int_{\Sph^{2}} \rho^{-1} \Psi^{j}_{lm}(\Omega) 
	\Psi^{j}_{l'm'}(\Omega)\, d^{2}\Omega = -\rho^{-1} \delta_{mm'}
	\delta_{ll'}.
\end{align*}
    Then
    the higher
$ l $
    terms with different momenta are still orthogonal
    in the product in the integral
(\ref{sterm})
    and they vanish due to the decrease of
$ u_{lm} $,
$ \phi_{lm} $
    at the origin, so we get
\begin{equation*}
    \int_{\partial B_{\rho}}
	\partial_{k} f^{j}(\vec{x}) f^{j}(\vec{x}) \,d^{2}\sigma^{k}
    \simeq 2\frac{u_{1m}u'_{1m}}{\rho^{2}} -4 \frac{u^{2}_{1m}}{\rho^{3}}
	- \frac{u'^{2}_{1m}}{\rho} + u'_{1m}u''_{1m} , \quad\rho\to 0,
\end{equation*}
    where
$ u_{1m} = u_{1m}(\rho) $.
    Together with the other terms in
(\ref{qftrans})
    this leads to the following expression
\begin{multline}
\label{qf2trans}
Q_{\Delta}^{\kappa}(\vec{f}) 
    = \lim_{\rho\to 0} \Bigl(
    \int_{\RR^{3}\setminus B_{\rho}}
	\bigl( \frac{\partial f^{j}(\vec{x})}{\partial x_{k}}
	\bigr)^{2} d^{3}x
    +2\frac{u_{1m}(\rho)u_{1m}'(\rho)}{\rho^{2}} -\\
    -4 \frac{u_{1m}^{2}(\rho)}{\rho^{3}} - \frac{u_{1m}'^{2}(\rho)}{\rho}
	+ u_{1m}'(\rho)u_{1m}''(\rho)
    - 2u_{1m}'(0) \frac{u_{1m}(\rho)}{\rho^{2}} \Bigr) .
\end{multline}
    Now we can write the Taylor expansions
\begin{equation*}
    u(\rho) = u'(0) \rho + u''(0)\frac{\rho^{2}}{2} + \ldots, \quad
    u'(\rho) = u''(0) \rho + \ldots, \quad
    \rho\to 0,
\end{equation*}
    where, according to
(\ref{bcond})
\begin{equation*}
    u''(0) = \frac{4}{3}\kappa u'(0),
\end{equation*}
    and substitute them into the boundary terms in the expression
(\ref{qf2trans})
    to obtain that
\begin{align*}
    \frac{u(\rho)u'(\rho)}{\rho^{2}} &\simeq \frac{1}{\rho}u'^{2}(0)
+ \frac{3}{2} u'(0)u''(0) = \bigl(\frac{1}{\rho} +2\kappa\bigr) u'^{2}(0) ,\\
    \frac{u^{2}(\rho)}{\rho^{3}} &\simeq \frac{u'^{2}(0)}{\rho} +u'(0)u''(0)
	= \bigl(\frac{1}{\rho}+\frac{4}{3}\kappa\bigr) u'^{2}(0) ,\\
    \frac{u'^{2}(\rho)}{\rho} &\simeq \frac{u'^{2}(0)}{\rho} +2u'(0)u''(0)
	= \bigl(\frac{1}{\rho}+\frac{8}{3}\kappa\bigr) u'^{2}(0) ,\\
    u'(\rho)u''(\rho) &\simeq \frac{4}{3}\kappa u'^{2}(0) ,\quad
    \frac{u(\rho)}{\rho^{2}} \simeq \frac{1}{\rho}u'(0)+\frac{1}{2}u''(0)
	= \bigl(\frac{1}{\rho} +\frac{2}{3}\kappa \bigr)u'(0)
\end{align*}
    as
$ \rho\to 0 $
    and hence
\begin{equation}
\label{uTu}
Q_{\Delta}^{\kappa}(\vec{f}) = \lim_{\rho \to 0} \Bigl(
    \int_{\RR^{3}\setminus B_{\rho}}
	\bigl( \frac{\partial f^{j}_{m}(\vec{x})}{\partial x_{k}}
	\bigr)^{2} d^{3}x
    - \bigl(\frac{5}{\rho} + 4\kappa \bigr)u'^{2}_{m}(0) \Bigr).
\end{equation}
    Similar expansion, written for the integral
\begin{equation*}
    \int_{\partial B_{\rho}} |\vec{f}|^{2}(\vec{x})\,d^{2}s
    \simeq 2\frac{u_{1m}^{2}(\rho)}{\rho^{2}} + u_{1m}'^{2}(\rho)
    \simeq \bigl(3 + \frac{16}{3}\rho \kappa\bigr) u_{1m}'^{2}(0) ,
	\quad \rho\to 0,
\end{equation*}
    can be substituted in the RHS of
(\ref{uTu})
    and this yields the following basis-independent equation
\begin{equation}
\label{mext}
Q_{\Delta}^{\kappa}(\vec{f}) 
    = \lim_{\rho\to 0}\Bigl(
    \int_{\RR^{3}\setminus B_{\rho}}
\bigl(\frac{\partial f^{j}(\vec{x})}{\partial x_{k}}\bigr)^{2} d^{3} x -
    \bigl(\frac{5}{3\rho}- \frac{44}{27}\kappa\bigr) \int_{\partial B_{\rho}}
        |\vec{f}|^{2}(\vec{x})\, d^{2} s \Bigr) .
\end{equation}
    It is worth to note, that this expression coincides, up to the definition
    of the parameter
$ \kappa $,
    with the quadratic form in
\cite{Lapl}.
    But now the form is extended and is closable \emph{w.~r.}~to the
    physical scalar product
(\ref{fgprod}).
    A non spherically symmetric extension, with different parameters
$ \kappa $
    at different
$ m $,
    requires some projectors separating the
$ l=1 $
    components of the function
$ \vec{f}(\vec{x}) $.

\section{The inverse operator $ T^{-1} $}
	To analyse the quadratic form
(\ref{t1sform})
	of the inverse operator
$ T^{-1} $
	let us express the latter form in terms of scalar product
(\ref{angleprod})
\begin{equation*}
    Q^{-1}(u) = (u,T^{-1}u) = (u,T T^{-2} u) = \langle u, T^{-2} u \rangle ,
\end{equation*}
	where
\begin{equation*}
    T^{-2}(r,s) = \int T^{-1}(r,t) T^{-1}(t,s) \, dt =
	\frac{1}{6}\bigl(r^{2}s -\frac{r^{4}}{5s}\bigr) \theta(s-r)
	+\frac{1}{6}\bigl(s^{2}r -\frac{s^{4}}{5r}\bigr) \theta(r-s) .
\end{equation*}
	We define a symmetric operator
$ \tilde{T}^{-2} = T^{-2}|_{\WW_{0}} $
	by means of the action of
$ T^{-2} $
	on a set
\begin{equation*}
    \WW_{0} = \{u(r):\; u\in\HH, \,
	\langle T^{-2}u, T^{-2}u \rangle < \infty, \,
	\text{+ conditions} \} ,
\end{equation*}
    where \emph{conditions} are
\begin{align}
\label{cr0}
    \int_{0}^{\infty} D^{*} u(r)\, dr &
	= - \int_{0}^{\infty} \frac{u(r)}{r} dr = 0 , \\
\label{cr1}
    \int_{0}^{\infty} r D^{*} u(r)\, dr & 
	= -r u(r)\bigr|_{0}^{\infty} = 0 , \\
\label{cr2}
    \int_{0}^{\infty} r^{2} D^{*} u(r)\, dr & 
	= \int_{0}^{\infty} r u(r) \,dr = 0 , \\
\label{cr3}
    \int_{0}^{\infty} r^{3} D^{*} u(r)\, dr & 
	= 2 \int_{0}^{\infty} r^{2} u(r) \,dr = 0 . 
\end{align}
	The set
$ \WW_{0} $
	is dense in
$ \HH $
	and closed in the graph norm of operator
$ \tilde{T}^{-2} $
	(the latter statement requires a separate justification).
	The symmetricity of
$ \tilde{T}^{-2} $
	follows from construction and smoothness
        of its kernel in
$ r $
	at fixed
$ s $,
\begin{equation*}
    \langle u, \tilde{T}^{-2} v \rangle = (u, TT^{-2}v) = (u,T^{-1}v)
	= (T^{-1}u, v) = \langle \tilde{T}^{-2} u, v\rangle .
\end{equation*}
	We also remark that the requirement
$ \langle T^{-2} u, T^{-2} u \rangle < \infty $
	yields the conditions
(\ref{cr1}) and
(\ref{cr3}),
	which will be relaxed when defining self-adjoint extenstions.

\subsection{Deficiency indices of operator
$ \tilde{T}^{-2} $}
	To calculate the deficiency indices let us use the following
	auxiliary formula
\begin{equation*}
    \int_{0}^{\infty} T^{-2}(r,s) D e^{\sigma s} \, ds
	= \frac{1}{\sigma^{4}} D(e^{\sigma r}-1) - \frac{r}{2\sigma^{2}},
    \quad \RE\,\sigma <0 .
\end{equation*}
	The following combinations of exponents falling off at infinity
	and vanishing at the origin can be used as vectors
$ c_{\pm} $
\begin{equation*}
    c_{\pm}(r) = \alpha_{\pm} D \exp\{e^{\mp\frac{5i}{8}\pi}\rho r\}
	+ \beta_{\pm} D \exp\{e^{\mp\frac{9i}{8}\pi}\rho r\}
	- \frac{\alpha_{\pm}+\beta_{\pm}}{r} 
\end{equation*}
    from the kernels of the adjoint operators
\begin{equation*}
    0 = \langle c_{\pm} , (\tilde{T}^{-2} \pm i\rho^{-4}) u \rangle =
    \bigl(Tc_{\pm} , (T^{-2}\pm i\rho^{-4}) u \bigr) ,\quad u\in\WW_{0} .
\end{equation*}
	Then using the above auxiliary formula this equation
	can be continued and checked as following
\begin{align*}
    \bigl(T &c_{\pm}, (T^{-2}\pm i\rho^{-4}) u \bigr) = \\
	=& \bigl( \rho^{2} (T^{-2}\mp i\rho^{-4}) 
(\alpha_{\pm} e^{\mp\frac{5i}{4}\pi} D \exp\{e^{\mp\frac{5i}{8}\pi}\rho r\}
+\beta_{\pm} e^{\mp\frac{i}{4}\pi} D \exp\{e^{\mp\frac{9i}{8}\pi}\rho r\} ) ,
    u \bigr) \\
    =& e^{\pm\frac{i\pi}{4}} \rho^{-2}
	(\alpha_{\pm}-\beta_{\pm}) \bigl(r^{-1}, u\bigr)
	+ \frac{1}{2} (\alpha_{\pm}+\beta_{pm}) \bigl(r, u\bigr) = 0,
\end{align*}
    where the last equality is derived from
(\ref{cr0}) and 
(\ref{cr2}).
    Thus the vectors
$ c_{\pm} $
    span 2-dimensional spaces and the deficiency indices of operator
$ \tilde{T}^{-2} $
    are (2,2).

\subsection{First example of self-adjoint extensions}
    A self-adjoint extension of symmetric operator with deficiency indices
    (2,2) is described by a point in the manifold of the unitary group
$ U(2) $.
    In what follows we restrict the presentation
    to the subgroup of this set, namely to
    the extensions defined on the domains
\begin{equation*}
    \WW_{\kappa} = \{u(r):\; u\in\HH,\, T^{-2}_{\kappa}u \in\HH, \,
\kappa^{3}\int_{0}^{\infty}r^{2}u\,dr = 6\int_{0}^{\infty}\frac{u}{r}\, dr \},
    \quad \kappa \in \RR .
\end{equation*}
    The condition in
$ \WW_{\kappa} $
    at nonzero
$ \kappa $
    eliminates the asymtotics
$ r^{-1} $
    at infinity, i.~e. it implies the condition
(\ref{cr1}).
	Similar to how resolvents (or the inverse operators) of
	self-adjoint extensions of a differential operator
	differ only by a finite-dimensional projector
(see ref.
\cite{AK}),
	two distinct self-adjoint extensions in question
	differ by a simple expression as well.
	Let us define the action of the operator
$ T^{-2}_{\kappa} $
	as the action of the kernel of
$ T^{-2} $
	augmented by a one-dimensional additive contribution,
\begin{equation*}
    T^{-2}_{\kappa}(r,s) = T^{-2}(r,s) - \frac{r}{\kappa^{3} s} .
\end{equation*}
	Symmetricity of 
$ T^{-2}_{\kappa} $
	follows from symmetricity of
$ T^{-2} $
	and that of
$ -\frac{r}{\kappa^{3}s} $
	with respect to product
$ \langle \cdot , \cdot \rangle $.
	Self-adjointness (given the proof of the structure of the kernels
of operator adjoint to $ \tilde{T}^{-2} $)
	now is inferred from the following relations which can be found
	using equations
(\ref{Pf2}), (\ref{Pf3})
\begin{equation*}
    \int_{0}^{\infty} \frac{d_{\pm}}{r} dr
	= \sqrt{2} e^{\pm\frac{5i}{8}\pi} \rho , \quad
    \int_{0}^{\infty} r^{2} d_{\pm} dr
	= 6 e^{\pm\frac{i}{4}\pi}\rho^{-2} ,
\end{equation*}
    where
\begin{equation*}
    d_{\pm}(r) = D \exp\{e^{\mp\frac{5i}{8}\pi}\rho r\}
	- D \exp\{e^{\mp\frac{9i}{8}\pi}\rho r\} .
\end{equation*}
    Vectors proportional to
$ d_{\pm} $
	do not obey the integral condition of
$ \WW_{\kappa} $
	at any real
$ \kappa $.

	Finally we would like to present a simple expression
	for the quadratic form of
$ T^{-2}_{\kappa} $,
	which does not involve scalar product
(\ref{plainprod})
\begin{equation*}
    Q^{-1}_{\kappa}(u) = \langle u , T^{-2}_{\kappa} u\rangle
	= (u,TT^{-2}_{\kappa} u) = \iint_{0}^{\infty} \ol{u(r)}\bigl(
    T^{-1}(r,s) - \frac{2}{\kappa^{3}rs}\bigr) u(s)\,dr\,ds .
\end{equation*}
	This expression shows that form
(\ref{t1sform})
	corresponds to the value
$ \kappa = \pm\infty $.
	At the same time, any form with a non-zero parameter
$ \kappa $
	is an extension of a ``maximal'' form (according to Friedrichs-Stone
theorem
\cite{FS})
$ Q^{-1}_{\kappa}|_{\kappa = -0} $
	defined on the set of functions that obey
(\ref{cr0}),
(\ref{cr1}).

\subsection{The resolvent of operator
$ T^{-2}_{\kappa} $}
	To examine the spectral properties of operator
$ T^{-2}_{\kappa} $
	we construct an expression for the kernel of its resolvent
$ R(r,s;z) $ in terms of variable
$ z = w^{-4} $:
\begin{equation}
\label{T2reseq}
    (T^{-2}_{\kappa} - z^{-4}) R(r,s;z) = \delta(r-s) ,
	\quad 0 < \arg z < \frac{\pi}{2} .
\end{equation}
	We demand that this kernel should satisfy the boundary conditions from
$ \WW_{\kappa} $
	in variable
$ r $,
	fall off at infinity in variables
$ r $,
$ s $
	and be complex symmetric in these variables with respect
	to scalar product
$ \langle \cdot , \cdot \rangle $
\begin{equation}
\label{CS}
    R^{*}(r,s;z) = R(s,r;i\ol{z}) .
\end{equation}
	We shall approach this task as follows.
	A formal equality
\begin{equation*}
    (T^{-2} - z^{-4})^{-1} = - \frac{z^{4}T^{2}}{T^{2}-z^{4}}
	= -\frac{z^{4}}{2}\bigl(\frac{1}{T-z^{2}}+\frac{1}{T+z^{2}}\bigr) T
\end{equation*}
	suggests that the resolvent can be sought as a sum
\begin{equation}
\label{Rprel}
    R(r,s;z) = -\frac{z^{4}}{2} \bigl(R_{-}(z)+R_{+}(z)+\tilde{R}(z)\bigr)T
\end{equation}
	where
$ R_{\pm} $ and
$ \tilde{R} $
	are symmetric functions such that
\begin{gather}
\label{RReqs1}
    (T\pm z^{2})R_{\pm} = \delta(r-s) ,\quad (T^{2}-z^{4}) \tilde{R} = 0 ,\\
\label{RReqs2}
    \ol{R_{\pm}(r,s;z)} = R_{\mp}(s,r;i\ol{z}) ,\quad
    \ol{\tilde{R}(r,s;z)} = \tilde{R}(s,r;i\ol{z}) .
\end{gather}
	In this representation, complex symmetricity 
(\ref{CS}) 
	is explicit from equation
(\ref{apc}).
	Let us pick expressions for
$ R_{\pm} $ and
$ \tilde{R} $
	in such a way that sum
(\ref{Rprel})
	would obey the boundary conditions in
$ r $.
	Functions
$ De^{\pm izr} $,
$ De^{\pm zr} $
	satisfy the homogeneous equations
\begin{equation*}
    (T + z^{2}) De^{\pm zr} = 0, \quad (T - z^{2}) De^{\pm izr} = 0 ,
\end{equation*}
	while at the same time
$ De^{izr} $,
$ De^{-zr} $
	fall off rapidly at
$ r\to\infty $
	in the sector
$ 0<\arg z <\frac{\pi}{2} $.
	We shall seek
$ R_{\pm} $ and
$ \tilde{R} $
	in the following form,
\begin{align}
\label{Rpm1}
    R_{+}(r,s;z) =&\, \frac{1}{W_{+}}\bigl(
	De^{zr}De^{-zs} \theta(s-r) + De^{zs} De^{-zr} \theta(r-s)
	\bigr),\\
\label{Rpm2}
    R_{-}(r,s;z) =&\, \frac{1}{W_{-}}\bigl(
	De^{-izr}De^{izs} \theta(s-r) + De^{-izs} De^{izr} \theta(r-s)
	\bigr),\\
\label{tR1}
    \tilde{R}(r,s;z) =&\, \frac{\beta_{-}(z)}{W_{-}} \bigl(
	De^{-zr}De^{izs} \theta(s-r) + De^{-zs}De^{izr}\theta(r-s) \bigr) \\
\label{tR2}
    &+ \frac{\beta_{+}(z)}{W_{+}} \bigl(
	De^{izr}De^{-zs} \theta(s-r) + De^{izs}De^{-zr}\theta(r-s) \bigr) \\
\label{tRa}
    &+ \frac{\alpha_{+}}{W_{+}}De^{-zr}De^{-zs}
	+ \frac{\alpha_{-}}{W_{-}}De^{izr}De^{izs},
\end{align}
    where
$ W_{\pm} $
    denote the Wronskians
\begin{align*}
    W_{+}(z) &= \frac{\partial}{\partial r} De^{zr} De^{-zr}
	- De^{zr} \frac{\partial}{\partial r} De^{-zr} = -2z^{3} ,\\
    W_{-}(z) &= \frac{\partial}{\partial r} De^{-izr} De^{izr}
	- De^{-izr} \frac{\partial}{\partial r} De^{izr} = -2iz^{3} .
\end{align*}
	Only the first terms in expressions
(\ref{Rpm1}), (\ref{Rpm2}), (\ref{tR1}), (\ref{tR2})
	are involved in the boundary conditions
	for resolvent
$ R(r,s;z) $
	in variable
$ r $
	at fixed
$ s $.
	Together with
(\ref{tRa})
	they can be grouped into two expressions
\begin{equation*}
    -\frac{z^{4}TDe^{-zs}}{2W_{+}}
	D(e^{zr} + \alpha_{+}e^{-zr} + \beta_{+} e^{izr}) 
    - \frac{z^{4}TDe^{izs}}{2W_{-}}
	D(e^{-izr} + \alpha_{-}e^{izr} + \beta_{-} e^{-zr}) ,
\end{equation*}
	which are proportional to different functions of
$ s $.
	The boundary conditions from
$ \WW_{\kappa} $
	must be fulfilled for either of these terms.
	Using expressions
(\ref{Pf2}), (\ref{Pf3}),
	this leads to the following equations for
$ \alpha_{\pm}(z) $,
$ \beta_{\pm}(z) $:
\begin{align*}
    1 + \alpha_{+} + \beta_{+} & = 1+ \alpha_{-} + \beta_{-} = 0 ,\\
\kappa^{3} (1+\alpha_{-}-\beta_{-})&= 2z^{3} (i - i\alpha_{-} + \beta_{-}), \\
\kappa^{3} (1+\alpha_{+}-\beta_{-})&= 2z^{3} (1 - \alpha_{+} + i\beta_{+}) .
\end{align*}
	These equations have a single solution
\begin{align}
\label{abm}
    \alpha_{-} &= - \frac{\kappa^{3}+z^{3}(1-i)}{\kappa^{3}+z^{3}(1+i)} ,
	\quad \beta_{-} = \frac{-2iz^{3}}{\kappa^{3}+z^{3}(1+i)}
	    = \frac{W_{-}(z)}{d(z)} ,\\
\label{abp}
    \alpha_{+} &= - \frac{\kappa^{3}+z^{3}(i-1)}{\kappa^{3}+z^{3}(1+i)} ,
	\quad \beta_{+} = \frac{-2z^{3}}{\kappa^{3}+z^{3}(1+i)}
	    = \frac{W_{+}(z)}{d(z)} ,
\end{align}
	where
$ d(z) $ is the common denominator (the determinant)
\begin{equation*}
    d(z) = \kappa^{3} + z^{3} (1+i) .
\end{equation*}
	The above explicit formulas for functions
$ \beta_{\pm}(z) $
	immediately lead us to a ``magic'' relation
\begin{equation*}
    \frac{\beta_{-}(z)}{W_{-}(z)} = \frac{\beta_{+}(z)}{W_{+}(z)}
	= \frac{1}{d(z)} ,
\end{equation*}
	which further allows to cancel the
$ \theta $-functions
	in the term
$ \tilde{R}(r,s;z) $
	and bring the latter to the following form
\begin{align*}
    \tilde{R}(r,s;z) =&\, \frac{1}{d(z)} \bigl(De^{-zr}De^{izs}
	+ De^{-zs} De^{izr} \bigr) \\
    &+ \frac{\alpha_{+}}{W_{+}}De^{-zr}De^{-zs}
	+ \frac{\alpha_{-}}{W_{-}}De^{izr}De^{izs}.
\end{align*}
	And this, in its turn, means that
$ R_{\pm} $ and
$ \tilde{R} $
	obey equations
(\ref{RReqs1})
	and relation
(\ref{RReqs2}).

	We, therefore, have shown that the conjectured solution
(\ref{Rprel})
	does satisfy the boundary conditions in
$ r $,
	its fall off at infinity following from the exponential
	decrease of functions
$ De^{izr} $,
$ De^{-zr} $.
	In order to ensure that
$ R(r,s;z) $
	indeed is the resolvent of operator
$ T^{-2}_{\kappa} $
	one has to check that it obeys equation
(\ref{T2reseq}).
	We do not show here the related calculations due to their high volume,
	but one can check that this equation is indeed fulfilled.
    
\subsection{Spectral properties of operator
$ T^{-2}_{\kappa} $}
	The identity decomposition
(\ref{Acompl}),
	written in terms of variable
$ z = w^{-1/4} $,
	which we have used in resolvent
$ R(r,s;z) $,
	turns into the following expression
\begin{align}
\label{T2compl}
    &\delta(r-s) = \frac{1}{2\pi i} \int_{0}^{\infty} \bigl(
	R(r,s;\lambda)-R(r,s;i\lambda) \bigr) \frac{4}{\lambda^{5}}\,
    d\lambda \\
\nonumber
    &- \sum_{n} \lim_{z\to z_{n}} \frac{4}{z_{n}^{5}} R(r,s;z) (z_{n}-z)
= \int_{0}^{\infty} \hat{p}_{\lambda}(r) T\hat{p}_{\lambda}(s) \,d\lambda
	+ \sum_{n} \hat{q}_{n}(r) T\hat{q}_{n}(s) ,
\end{align}
	while the orthogonality relations preserve their form 
(\ref{ortrel1}), (\ref{ortrel2}).
	Resolvent
$ R(r,s;z) $
	can only have one pole --- it is determined by the zero of the denominator
$ d(z) $
	and located at the point
$ z_{0} = 2^{-1/6} e^{i\pi/4} \kappa $.
	In order for this pole to fall into the sector
$ 0 < \arg z < \frac{\pi}{2} $
	it is necessary that
$ \kappa > 0 $.
	If that is the case, given that the terms
$ R_{\pm}(z) $
	are regular in
$ z $,
	the residue of the resolvent can be calculated and transformed as follows,
\begin{align*}
    \frac{4}{z_{0}^{5}}& \lim_{z\to z_{0}} R(r,s;z)(z-z_{0}) =
	-\frac{2}{z_{0}} \lim_{z\to z_{0}} \tilde{R}(r,s;z)(z-z_{0}) T_{s} = \\
    =&\, -\frac{2}{z_{0}}\lim_{z\to z_{0}} \bigl( \frac{z-z_{0}}{d}
	(De^{-zr}De^{izs}+De^{-zs}De^{izr})
	+\frac{\alpha_{+}(z-z_{0})}{W_{+}}De^{-zr}De^{-zs} \\
	&+\frac{\alpha_{-}(z-z_{0})}{W_{-}}De^{izr}De^{izs} \bigr) T_{s}
    = \frac{-2 T_{s}}{3(i+1)z_{0}^{3}}
	(De^{-z_{0}r}De^{iz_{0}s}+De^{-z_{0}s}De^{iz_{0}r}) \\
    &- \frac{4}{3i(i+1)^{3}z_{0}^{3}}
	(De^{-z_{0}r}De^{-z_{0}s}+De^{iz_{0}r}De^{iz_{0}s}) \\
    =& \frac{-2}{3\kappa^{3}}
    \bigl(D\exp\{2^{-\frac{1}{6}}e^{\frac{3i}{4}\pi}\kappa r\}
	-D\exp\{2^{-\frac{1}{6}}e^{-\frac{3i}{4}\pi}\kappa r\} \bigr) 
    T_{s} \bigl(D\exp\{2^{-\frac{1}{6}}e^{\frac{3i}{4}\pi}\kappa s\} \\
	&-D\exp\{2^{-\frac{1}{6}}e^{-\frac{3i}{4}\pi}\kappa s\} \bigr) .
\end{align*}
	The latter product corresponds to the following (real) eigenfunction
	of the simple discrete spectrum
\begin{equation}
\label{T2qr}
    \hat{q}(r) = i\sqrt{\frac{2}{3\kappa^{3}}}
	\bigl(D\exp\{2^{-\frac{1}{6}}e^{\frac{3i}{4}\pi}\kappa r\}
	-D\exp\{2^{-\frac{1}{6}}e^{-\frac{3i}{4}\pi}\kappa r\} \bigr) .
\end{equation}

	The integral in the first line of the decomposition
(\ref{T2compl})
	corresponds to the continuous component of the spectrum.
	Let us first expand those parts of the resolvent that contain
$ R_{\pm}(z) $,
\begin{align}
\label{T2Rpm}
    R_{+}(&\, r,s;\lambda) + R_{-}(r,s;\lambda)
    - R_{+}(r,s;i\lambda) - R_{-}(r,s;i\lambda) = \\
\nonumber
    =&\, \frac{1}{2\lambda^{3}} \bigl(
	-De^{\lambda r}De^{-\lambda s} \theta(s-r)
	-De^{\lambda s}De^{-\lambda r} \theta(r-s) 
	+iDe^{-i\lambda r}De^{i\lambda s} \theta(s-r) \\
\nonumber
    &+iDe^{-i\lambda s}De^{i\lambda r} \theta(r-s) 
	+iDe^{i\lambda r}De^{-i\lambda s} \theta(s-r) 
	+iDe^{i\lambda s}De^{-i\lambda r} \theta(r-s) \\
\nonumber
    &+De^{\lambda r}De^{-\lambda s} \theta(s-r)
	+De^{\lambda s}De^{-\lambda r} \theta(r-s) \bigr) \\
\nonumber
    =&\, \frac{i}{2\lambda^{3}} \bigl(
	De^{-i\lambda r}De^{i\lambda s} + De^{i\lambda r}De^{-i\lambda s}
    \bigr).
\end{align}
	In the last transformation here we have cancelled the first two terms with the last two terms,
	while in the rest of the terms we have changed the sums of 
$ \theta $-functions
	of opposite arguments to unity.
	To transform the term
$ \tilde{R}(r,s;z) $
	we use the following relations for the solutions
(\ref{abm}),
(\ref{abp})
\begin{gather*}
    \ol{d(i\lambda)} = d(\lambda) \equiv d, \quad
	\alpha_{-}(\lambda) = - \frac{\ol{d}}{d} , \quad
	\alpha_{+}(i\lambda) = - \frac{d}{\ol{d}} , \\
\alpha_{-}(i\lambda) - \alpha_{+}(\lambda) = \frac{4i\lambda^{6}}{d\ol{d}}.
\end{gather*}
	Then
\begin{align*}
    \tilde{R}(r&,s;\lambda) - \tilde{R}(r,s;i\lambda) =\\
    =&\,\frac{1}{2\lambda^{3}} \bigl(
	\frac{2\lambda^{3}}{d}
	    (De^{-\lambda r}De^{i\lambda s} + De^{-\lambda s}De^{i\lambda r})
    -\alpha_{+}(\lambda) De^{-\lambda r} De^{-\lambda s} \\
    &+i\alpha_{-}(\lambda) De^{i\lambda r} De^{i\lambda s}
	- \frac{2\lambda^{3}}{\ol{d}}
	(De^{-i\lambda r}De^{-\lambda s} + De^{-i\lambda s}De^{-\lambda r})\\
    &+i\alpha_{+}(i\lambda) De^{-i\lambda r} De^{-i\lambda s}
    +\alpha_{-}(i\lambda) De^{-\lambda r} De^{-\lambda s} \bigr) \\
    =&\,\frac{i}{2\lambda^{3}} \bigl(
	\frac{2i\lambda^{3}}{\ol{d}}
	(De^{-i\lambda r}De^{-\lambda s} + De^{-i\lambda s}De^{-\lambda r})
	    -\frac{\ol{d}}{d} De^{i\lambda r} De^{i\lambda s} \\
    &-\frac{2i\lambda^{3}}{d}
	    (De^{-\lambda r}De^{i\lambda s} + De^{-\lambda s}De^{i\lambda r})
	    -\frac{d}{\ol{d}} De^{-i\lambda r} De^{-i\lambda s}
    +\frac{4\lambda^{6}}{d\ol{d}} De^{-\lambda r} De^{-\lambda s} \bigr) .
\end{align*}
	Taking into account the contribution
(\ref{T2Rpm})
	and the coefficient
$ 4T_{s}/(2\pi i \lambda^{5}) $
	these terms factorize into two instances of a single function
        taken at different points
$ r $
	and
$ s $
\begin{align*}
    \frac{4}{2\pi i\lambda^{5}} \bigl(R(r,s;\lambda) &-R(r,s;i\lambda)\bigr) 
    = \frac{-1}{2\pi\lambda^{4}}
	D\Bigl(\sqrt{\frac{\ol{d}}{d}}e^{i\lambda r}
	    - \sqrt{\frac{d}{\ol{d}}}e^{-i\lambda r}
	    +\frac{2i\lambda^{3}}{|d|}e^{-\lambda r}\Bigr) \\
	&\times T D\Bigl(\sqrt{\frac{\ol{d}}{d}}e^{i\lambda s}
	    - \sqrt{\frac{d}{\ol{d}}}e^{-i\lambda s}
	    +\frac{2i\lambda^{3}}{|d|}e^{-\lambda s}\Bigr) .
\end{align*}
	We, therefore, obtain the following real expression
	for the kernel
$ \hat{p}_\lambda(r) $
	from expansion
(\ref{T2compl}),
\begin{equation}
\label{T2pl}
    \hat{p}_\lambda(r) = \frac{i}{\sqrt{2\pi}\lambda^{2}}
	D\Bigl(\sqrt{\frac{\ol{d}}{d}}e^{i\lambda r}
	    - \sqrt{\frac{d}{\ol{d}}}e^{-i\lambda r}
	    +\frac{2i\lambda^{3}}{|d|}e^{-\lambda r}\Bigr) .
\end{equation}

\subsection{Second example of self-adjoint extensions}
    The second notable example of self-adjoint extensions of operator
$ \tilde{T}^{-2} $
    is based on the decompositions and eigensets
(\ref{Tqr}) and
(\ref{Tpl})
    of the previous section.
    The domains in question are
\begin{align*}
    \WW_{\varkappa} = \{u(r):\; u\in\HH,\, T^{-2}_{\varkappa}u \in\HH, & \,
\int_{0}^{\infty} D^{*}u\,dr = -\varkappa \int_{0}^{\infty} rD^{*}u\,dr , \\
	& 3\int_{0}^{\infty} r^{2} D^{*}u\,dr
	    = -\varkappa \int_{0}^{\infty} r^{3} D^{*}u\, dr \} ,
\end{align*}
    and the extended operators are defined as the following mixed kernels
\begin{equation*}
    T_{\varkappa}^{-2} = 
	\frac{1}{6}\bigl(s^{2}r -\frac{s^{4}}{5r}\bigr) \theta(r-s)
	+\frac{r}{12} s^{2} \frac{d}{ds} s - \frac{r^{2}}{6} s\frac{d}{ds}s
	+ \frac{r^{3}}{8} \frac{d}{ds}s .
\end{equation*}
    Corresponding quadratic forms can be written as one symmetric expression
\begin{equation*}
    Q_{\varkappa}^{-1}(u) =
	\langle (T^{-1}(r,s) - \frac{r}{2} \frac{d}{ds}s)u(s) ,
	(T^{-1}(r,s') - \frac{r}{2}\frac{d}{ds'}s')u(s') \rangle ,
\end{equation*}
    with 
$ s $ and $ s' $
    being integrated.
    It does not explicitly contain the parameter
$ \varkappa $,
    but the latter enters the domains of the forms as the condition
\begin{equation*}
\int_{0}^{\infty} D^{*}u\,dr = -\varkappa \int_{0}^{\infty} rD^{*}u\,dr .
\end{equation*}
    The eigensets
(\ref{Tqr}) and
(\ref{Tpl})
    allows to write the following spectral decomposition
\begin{equation*}
    Q_{\varkappa}^{-1}(u) = \iint_{0}^{\infty} Q_{\varkappa}^{-1}(r,s) u(r)
	\overline{u(s)} \,dr\,ds ,
\end{equation*}
    where
\begin{equation*}
    Q_{\varkappa}^{-1}(r,s) = \int_{0}^{\infty} T_{r}
	\tilde{p}_{\lambda,\varkappa}(r)
    T_{s}\tilde{p}_{\lambda,\varkappa}(s) \lambda^{-4} \,d\lambda
	- \varkappa^{-4} T_{r}\tilde{q}_{\varkappa}(r)
	    T_{s}\tilde{q}_{\varkappa}(s) .
\end{equation*}

\section{Conclusion}
	We have studied the spectral properties of self-adjoint extensions
	of the $ l=1 $ radial part of the Laplace operator
	and those of the inverse operator to the latter part,
	in scalar product
(\ref{angleprod}).
	In both cases there is a simple continuous spectrum that
        occupies the non-negative half-axis and is described by
	the spectral densities (``eigenfunctions of the continuous spectrum'')
        presented in equations
(\ref{Tpl}) and
(\ref{T2pl}).
	The corresponding identity decompositions are described by expressions
(\ref{Tcompl}) and
(\ref{T2compl}).
	At negative (or positive, in the case of the inverse operator)
	values of the extension parameter
$ \kappa $ 
	the operators possess simple eigenvalues
$ -\kappa^{2} $
($ -2^{2/3}\kappa^{-4} $)
	and eigenfunctions
(\ref{Tqr}) and
(\ref{T2qr}), correspondingly.

\section*{Acknowledgments}
The work is partially supported by the grant 14-11-00598
of Russian Science Foundation.

\end{document}